\documentclass[oneside,reqno,A4]{amsart}

\usepackage{latexsym}
\usepackage{amsmath,amsthm,amssymb,amsxtra}

\newcommand{\qdn}{\hspace*{-1.5mm}}
\newcommand{\qqdn}{\hspace*{-2.5mm}}

\newcommand{\+}{&\qqdn}%

\newcommand{\ffnk}[4]{\left[\qdn\ba{#1}#3\\#4\ea{\!\Big|\:#2}\right]}

\newcommand{\zbnm}[2]{\genfrac{(}{)}{0mm}{2}{#1}{#2}}

\newcommand{\sbnq}[2]{\Bigl[\!\ba{c}\!#1\!\\#2\ea\!\Bigr]}

\newcommand{\be}{\begin{equation}}
\newcommand{\ee}{\end{equation}}
\newcommand{\ba}{\begin{array}}
\newcommand{\ea}{\end{array}}
\newcommand{\bmn}{\begin{eqnarray}}
\newcommand{\emn}{\end{eqnarray}}
\newcommand{\bnm}{\begin{eqnarray*}}
\newcommand{\enm}{\end{eqnarray*}}
\newcommand{\bln}{\begin{subequations}}
\newcommand{\eln}{\end{subequations}}

\newcommand{\eqn}[1]{\begin{equation}#1\end{equation}}

\newcommand{\alp}{\alpha}
\newcommand{\bet}{\beta}
\newcommand{\gam}{\gamma}

\newtheorem{thm}{Theorem}

\newcommand{\thank}[2]{\begin{center}\parbox{#1}
{{\sc\bf Acknowledgement}:\,{\small\it#2}}\end{center}}

\newcommand{\bbtm}[4]{\bibitem{kn:#1}{#2,}~\emph{#3,}~{#4.}}	
\newcommand{\cito}[1]{\cite{kn:#1}}	
\newcommand{\citu}[2]{\cite[#2]{kn:#1}}

\begin{document} \vspace*{-15mm}

\title{Finite Form of the Quintuple Product Identity}
\author{William Y. C. Chen - Wenchang Chu - Nancy S. S. Gu}
\dedicatory{Center for Combinatorics: LPMC\\
	Nankai University, Tianjin 300071\\
        People's Republic of China\\[-3mm]}
\thanks{\qquad
Email addresses: \emph{chen@nankai.edu.cn},
		 \emph{chu.wenchang@unile.it}
	     and \emph{gu@nankai.edu.cn}}
\date{}

\begin{abstract}
The celebrated quintuple product identity 
follows surprisingly from an almost-trivial 
algebraic identity, which is the limiting 
case of the terminating $q$-Dixon formula.  
\end{abstract}

\maketitle\thispagestyle{empty}
\markright{Chen - Chu - Gu: Finite Form of the Quintuple Product Identity}  
\renewcommand{\thethm}{}\vspace{-9mm}


The celebrated quintuple product identity discovered by 
Watson~\cito{watson} (cf.\:\citu{gasper}{P\:147} also) 
states that 
\eqn{\sum_{k=-\infty}^{+\infty}
(1-xq^k)q^{3\zbnm{k}{2}}(qx^3)^{k}
\:=\:\label{quint}
[q,x,q/x;q]_{\infty}
[qx^2,q/x^2;q^2]_{\infty}
\quad\text{for}\quad|q|<1} 
where the $q$-shifted factorial is defined by
\[(x;q)_0\:=\:1
\quad\text{and}\quad
(x;q)_n=(1-x)(1-qx)\cdots(1-q^{n-1}x)
\quad\text{for}\quad
n=1,2,\cdots\]
with the following abbreviated multiple parameter notation
\[[\alp,\bet,\cdots,\gam;q]_\infty
\:=\:(\alp;q)_\infty
(\bet;q)_\infty
\cdots(\gam;q)_\infty.\]
This identity has several important applications in combinatorial 
analysis, number theory and special functions. For the historical 
note, we refer the reader to the paper~\cito{carlitz}.
In this short note, we shall show that identity \eqref{quint} 
follows surprisingly from the following algebraic identity.
\begin{thm}[Finite form of the quintuple product identity]
For a natural number $m$ and a variable $x$, there holds
an algebraic identity: 
\eqn{1\:\equiv\:\sum_{k=0}^m(1+xq^k)\sbnq{m}{k}
\frac{(x;q)_{m+1}}{(q^kx^2;q)_{m+1}}
x^kq^{k^2}.\label{quint-f}}
\end{thm}\vspace{-5mm}
In fact, performing parameter replacements $m\to m+n$, 
$x\to-q^{-m}x$ and $k\to k+m$ and then simplifying 
the result through factorial-fraction relation 
\bnm
&&\frac{(-q^{-m}x;q)_{m+n+1}}
     {(q^{k-m}x^2;q)_{m+n+1}}
\:=\:
\frac{(-q^{-m}x;q)_m(-x;q)_{1+n}}
     {(q^{k-m}x^2;q)_{m-k}(x^2;q)_{1+n+k}}\\
&=&(-1)^{m-k}q^{\zbnm{k}{2}-mk}x^{2k-m}
\times
\frac{(-q/x;q)_m(-x;q)_{1+n}}
     {(q/x^2;q)_{m-k}(x^2;q)_{1+n+k}}
\enm
we may restate the algebraic identity displayed in the 
theorem as the finite bilateral series identity
\eqn{1\:\equiv\:\sum_{k=-m}^n(1-xq^k)\sbnq{m+n}{m+k}
\frac{(-x;q)_{1+n}(-q/x;q)_{m}}
     {(x^2;q)_{1+n+k}(q/x^2;q)_{m-k}}
x^{3k}q^{k^2+\zbnm{k}{2}}.}
Letting $m,n\to\infty$ in this equation and applying
the relation 
\[(q;q)_{\infty}
\frac{(x^2;q)_{\infty}(q/x^2;q)_{\infty}}
     {(-x;q)_{\infty}(-q/x;q)_{\infty}}
\:=\:[q,x,q/x;q]_{\infty}
[qx^2,q/x^2;q^2]_{\infty}\]
we derive immediately the quintuple product identity 
displayed in \eqref{quint}.

In terms of basic hypergeometric series, we remark that 
the finite sum identity \eqref{quint-f} is just 
the limiting case $M\to\infty$ of the terminating 
$q$-Dixon formula (cf.\:\citu{gasper}{II-14}):
\[{_4\phi_3}
\ffnk{crccrc}{q;\frac{q^{1+m}x}{M}}
     {x^2,\+-qx,\+q^{-m},\+M}
     {\+-x,\+q^{1+m}x^2,\+qx^2/M}
\:=\:\frac{(qx^2;q)_m(qx/M;q)_m}
          {(qx;q)_m(qx^2/M;q)_m}.\]

\thank{150mm}{ This work was done under the 
auspices of the ``973'' Project on Mathematical Mechanization, 
the National Science Foundation, the Ministry of Education, 
and the Ministry of Science and Technology of China.}
\end{document}